\documentclass{svproc}

\usepackage{url}

\usepackage{amsmath}
\usepackage{amssymb}
\usepackage{booktabs}

\newcommand{\norm}[1]{\left\lVert #1 \right\rVert}

\begin{document}
\mainmatter

%======================================================================
% Title and author information
%======================================================================

\title{A Reduced Coupling-Space Multigrid Schur Preconditioner
for Fictitious-Domain Interface Problems}

\titlerunning{A Reduced Coupling-Space Multigrid Schur Preconditioner}

\author{Najwa Alshehri}

\authorrunning{N. Alshehri}

\tocauthor{Najwa Alshehri}

\institute{ Jubail Industrial College, 
Royal Commission for Jubail and Yanbu, Saudi Arabia\\
\email{shehrin@rcjy.edu.sa}
}
\maketitle
% 
%======================================================================
% Abstract
%======================================================================
% 
\begin{abstract}
We study a reduced multigrid Schur preconditioner for the fictitious-domain
method with distributed Lagrange multipliers applied to elliptic interface
problems. The method replaces the costly background solve by a multigrid
approximation and treats the resulting correction only in the multiplier
coupling space through a Woodbury formula. Numerical experiments show a marked
improvement over the baseline preconditioner and consistently rapid
convergence under mesh refinement and coefficient contrast.

\end{abstract}
% 
%======================================================================
\section{Introduction}
%======================================================================
Elliptic interface problems describe a physical field governed by elliptic partial differential equations in two or more subdomains separated by an internal interface, across which material coefficients may change and the solution must satisfy appropriate transmission conditions. They are important because they model phenomena such as heat conduction in composite materials, diffusion in heterogeneous media, and deformation in layered structures. They also provide a useful first step toward more complicated multiphysics problems, since they retain the main difficulty of enforcing conditions across an interface while avoiding time dependence, nonlinear effects, and moving domains.

To solve such problems without fitting the background mesh to the interface, we use the unfitted fictitious-domain method with distributed Lagrange multipliers (FD--DLM), introduced in a finite-element immersed-boundary framework in \cite{boffi2015finite}. Its formulation for elliptic interface problems was studied in \cite{boffi2014mixed,auricchio2015fictitious,Najwa2022elliptic}, where stable finite-element discretizations with continuous and discontinuous multiplier spaces were developed. A main benefit of the FD--DLM approach is that the meshes can remain fixed and independent of the interface, avoiding mesh regeneration, especially in time-dependent problems with moving interfaces.

The FD--DLM discretization leads to a large indefinite saddle-point system, so an effective preconditioner is essential. Earlier work considered block-diagonal and block-triangular preconditioners with exact block inversions \cite{boffi2024parallel}. The block-diagonal choice neglects important coupling and can give weak convergence. Multigrid was later used to improve the block inversions and gave better results than the diagonal approach \cite{alshehri2025multigrid}, although further improvement is still needed. A different direction is provided by augmented-Lagrangian preconditioners \cite{benzi2025optimal,benzi2026augmented}. In ongoing work with Boffi and Credali, we investigate a transforming smoother that recasts the system into a form with a more favorable right preconditioner. However, this still leaves the challenge of approximating the Schur complement accurately at low cost. The question addressed in this paper is therefore whether the background contribution to the Schur complement can be captured by a multigrid approximation in a reduced coupling space.

The remainder of the paper is organized as follows. Section~2 presents the elliptic interface problem, its FD--DLM formulation, and the discrete spaces. Section~3 introduces the reduced multigrid Schur preconditioner. Section~4 reports the numerical results.
% 
%======================================================================
\section{Problem Setting and FD--DLM Formulation}
%======================================================================
Let $\Omega\subset\mathbb{R}^d$, $d\in\{2,3\}$, be a bounded Lipschitz
domain with boundary $\partial\Omega$. We assume that $\Omega$ is divided
into two subdomains $\Omega_1$ and $\Omega_2$, separated by a Lipschitz
interface $\Gamma=\partial\Omega_1\cap\partial\Omega_2$. For simplicity,
$\Omega_2$ is fully immersed in $\Omega$, so that
$\Gamma\cap\partial\Omega=\varnothing$.

We consider an elliptic interface problem in two materials occupying
$\Omega_1$ and $\Omega_2$. Let $f_i\in L^2(\Omega_i)$ and let the coefficients
$\beta_i\in L^\infty(\Omega_i)$ satisfy
$0<\beta_{\min}\leq\beta_i(x)$ for a.e.
$x\in\Omega_i$, $i=1,2$. The unknowns $u_1$ and $u_2$ describe the solution
in the background and immersed subdomains, respectively. The problem is to
find $(u_1,u_2)$ such that
\begin{align}
 -\nabla\!\cdot(\beta_i\nabla u_i)&=f_i
 &&\text{in }\Omega_i, \qquad i=1,2, \label{eq:strong-diffusion}\\
 u_1&=u_2
 &&\text{on }\Gamma, \label{eq:strong-continuity}\\
 \beta_1\nabla u_1\cdot n_1+
 \beta_2\nabla u_2\cdot n_2&=0
 &&\text{on }\Gamma, \label{eq:strong-flux}\\
 u_1&=0
 &&\text{on }\partial\Omega. \label{eq:strong-bc}
\end{align}
Here, $n_i$ denotes the unit normal pointing outward from $\Omega_i$ on the
interface. Equation~\eqref{eq:strong-diffusion} gives the elliptic equation in
each subdomain. The first transmission condition,
\eqref{eq:strong-continuity}, enforces continuity of the solution across the
interface, while \eqref{eq:strong-flux} balances the normal fluxes. Finally,
\eqref{eq:strong-bc} prescribes a homogeneous Dirichlet condition on the outer
boundary of the background domain.

To obtain the FD--DLM formulation, we extend $f_1$ and $\beta_1$ smoothly from
$\Omega_1$ to the fictitious part of the background domain, and denote the
extensions by $f$ and $\beta$. The field $u$ represents the extension of
$u_1$ to $\Omega$; on the immersed domain $\Omega_2$, this extension is
required to coincide with $u_2$.

For a subdomain $\omega$, let $L^2(\omega)$ and $H^1(\omega)$ denote the
standard Lebesgue and Sobolev spaces. We set
$V=H^1_0(\Omega)$, where functions have zero trace on $\partial\Omega$,
$V_2=H^1(\Omega_2)$, and $\Lambda=V_2'$. The constraint
$u|_{\Omega_2}=u_2$ is imposed weakly through a distributed Lagrange
multiplier $\lambda\in\Lambda$. The continuous FD--DLM problem is to find
$(u,u_2,\lambda)\in V\times V_2\times\Lambda$ such that
\begin{align}
 (\beta\nabla u,\nabla v)_\Omega
 +\langle\lambda,v|_{\Omega_2}\rangle_{V_2',V_2}
 &= (f,v)_\Omega
 &&\forall v\in V, \label{eq:fddlm-background}\\
 ((\beta_2-\beta)\nabla u_2,\nabla v_2)_{\Omega_2}
 -\langle\lambda,v_2\rangle_{V_2',V_2}
 &= (f_2-f,v_2)_{\Omega_2}
 &&\forall v_2\in V_2, \label{eq:fddlm-immersed}\\
 \langle\mu,u|_{\Omega_2}-u_2\rangle_{V_2',V_2}&=0
 &&\forall\mu\in\Lambda. \label{eq:fddlm-constraint}
\end{align}
In other words, the multiplier tests the difference between the restriction of
the background field and the immersed field through the natural duality between
$H^1(\Omega_2)$ and its dual $[H^1(\Omega_2)]'$. Other choices of multiplier
space and coupling, including
an $H^1(\Omega_2)$ inner-product coupling for continuous multipliers, have
also been studied in \cite{auricchio2015fictitious,boffi2014mixed}.

We now pass from the continuous formulation to its finite-element
discretization. Let $\mathcal{T}_h$ and $\mathcal{T}_{2,h}$ be two independent
shape-regular quadrilateral meshes of $\Omega$ and $\Omega_2$, respectively,
with mesh sizes
$h=\max_{K\in\mathcal{T}_h}\operatorname{diam}(K)$ and
$h_2=\max_{K\in\mathcal{T}_{2,h}}\operatorname{diam}(K)$.
No matching condition is imposed between these meshes: elements of
$\mathcal{T}_h$ may be cut by $\Gamma$, whereas
$\mathcal{T}_{2,h}$ resolves the immersed domain independently.

We use the $Q_1$--$(Q_1+B)$--$P_0$ discretization. The background field is
approximated by continuous bilinear functions that vanish on
$\partial\Omega$. For the immersed field, we enrich the bilinear space with
bubble functions. The multiplier is approximated by discontinuous piecewise
constants. More precisely,
\[
V_h=Q_1(\mathcal{T}_h)\cap H^1_0(\Omega), \qquad
V_{2,h}=Q_1(\mathcal{T}_{2,h})\oplus B(\mathcal{T}_{2,h}), \qquad
\Lambda_h=P_0(\mathcal{T}_{2,h}).
\]
Here, $B(\mathcal{T}_{2,h})$ is the bubble-function space, whose functions
vanish on element boundaries. Since $\Lambda_h\subset L^2(\Omega_2)$, we use
the $L^2(\Omega_2)$ scalar product for the discrete coupling. We therefore
seek $(u_h,u_{2,h},\lambda_h)\in V_h\times V_{2,h}\times\Lambda_h$ such that
\begin{align}
 (\beta\nabla u_h,\nabla v_h)_\Omega
 +(\lambda_h,v_h|_{\Omega_2})_{\Omega_2}
 &= (f,v_h)_\Omega
 &&\forall v_h\in V_h, \label{eq:discrete-background}\\
 ((\beta_2-\beta)\nabla u_{2,h},\nabla v_{2,h})_{\Omega_2}
 -(\lambda_h,v_{2,h})_{\Omega_2}
 &= (f_2-f,v_{2,h})_{\Omega_2}
 &&\forall v_{2,h}\in V_{2,h}, \label{eq:discrete-immersed}\\
 (\mu_h,u_h|_{\Omega_2}-u_{2,h})_{\Omega_2}&=0
 &&\forall\mu_h\in\Lambda_h. \label{eq:discrete-constraint}
\end{align}
For the coefficient regime $\beta_2>\beta>0$, this stable element has been
shown to yield a unique discrete solution; see \cite{Najwa2022elliptic}.
% 
%======================================================================
\section{Preconditioner}
%======================================================================
The $Q_1$--$(Q_1+B)$--$P_0$ discretization considered above deserves particular
attention. Its discontinuous multiplier space leads to the
$L^2(\Omega_2)$ coupling, which is the more demanding coupling choice from the
conditioning viewpoint. For a related FD--DLM fluid--structure formulation,
the condition number associated with this coupling satisfies
$\operatorname{cond}(\mathcal A)\lesssim h^{-4}$ for comparable mesh sizes
under uniform refinement \cite{boffi2026stability}. This result motivates the
present study of the same coupling choice for elliptic interface problems and
provides a stringent test case for the proposed preconditioner.

After ordering the background degrees of freedom first, the discrete system has
the block form
\[
\begin{bmatrix}
A_1 & C_1^{T}\\
C_1 & B
\end{bmatrix}
\begin{bmatrix}
x_1\\
x_2
\end{bmatrix}
=
\begin{bmatrix}
f_1\\
f_2
\end{bmatrix}.
\]
Here, $x_1$ collects the degrees of freedom associated with $u_h$, whereas
$x_2$ collects those associated with $(u_{2,h},\lambda_h)$. The block $A_1$ is
the background elliptic block, $B$ contains the immersed-field and multiplier
variables, and $C_1$ describes their coupling.

Eliminating $x_1$ gives the Schur complement
$S=B-C_1A_1^{-1}C_1^{T}$. The exact Schur complement is useful: in ongoing
work, its use in a right Schur preconditioner yields mesh-independent
conditioning for the outer system. Its application, however, requires an
inner solve with $S$. Preconditioning that inner problem by $B^{-1}$ does not
give mesh-independent behavior, since the condition number of $SB^{-1}$
still increases under mesh refinement. This motivates the search for an
accurate, lower-cost approximation of $S^{-1}$, together with a cheaper
replacement for $A_1^{-1}$.

We take a different route. Instead of applying the exact background inverse,
which requires a global solve with $A_1$, we replace $A_1^{-1}$ by one
multigrid V-cycle $M_{\rm MG}$ and define
$S_{\rm MG}=B-C_1M_{\rm MG}C_1^T$. Although the correction
$C_1M_{\rm MG}C_1^T$ can be dense, it does not act on the whole space.
Indeed, writing $x_2=[x_{u_2};x_\lambda]$, the coupling has the form
$C_1^T=[\,0\;\;C_{1,\lambda}^T\,]$, and therefore
\[
S_{\rm MG}
=
B-
\begin{bmatrix}
0 & 0\\
0 & C_{1,\lambda}M_{\rm MG}C_{1,\lambda}^T
\end{bmatrix}.
\]
Thus, only the multiplier--multiplier block is changed. The background
response is consequently approximated only in the active multiplier coupling
space.

Let $r=\operatorname{rank}(C_{1,\lambda})$, and let
$Q_\lambda\in\mathbb{R}^{n_\lambda\times r}$ have orthonormal columns spanning
$\operatorname{range}(C_{1,\lambda})$. We embed this basis in the full Schur
space by setting
\[
Q=
\begin{bmatrix}
0\\
Q_\lambda
\end{bmatrix}.
\]
Thus, $Q$ represents precisely the multiplier directions that are coupled to
the background field. Since the multigrid correction has range contained in
this space, it can be written as
\[
C_1M_{\rm MG}C_1^T=QK_{\rm MG}Q^T, \qquad
K_{\rm MG}=Q^TC_1M_{\rm MG}C_1^TQ.
\]
Hence, $
S_{\rm MG}=B-QK_{\rm MG}Q^T$.
The proposed preconditioner applies its inverse through the reduced Woodbury
formula
 $
S_{\rm MG}^{-1}
=
B^{-1}+B^{-1}Q
\bigl(K_{\rm MG}^{-1}-Q^TB^{-1}Q\bigr)^{-1}Q^TB^{-1}$.
The only dense inverse is therefore of size $r\times r$, rather than the full
Schur-system size. The method captures the background response without forming
the dense Schur correction or applying the exact background inverse. In the
experiments below, $B^{-1}$ is applied exactly in order to isolate the effect
of the proposed correction; constructing a scalable approximation of $B^{-1}$
is left for future work.
% 
%======================================================================
\section{Numerical Results}
%======================================================================
We consider a two-dimensional test with background domain
$\Omega=[-1.4,1.4]^2$ and immersed domain $\Omega_2=B(0,1)$. The source
terms are $f=f_2=1$, and the background coefficient is $\beta=1$. We test the
two contrasts $\beta_2=10$ and $\beta_2=100$. The independent background and
immersed meshes are uniformly refined through four levels, with
$h=0.8879$, $0.4522$, $0.2278$, and $0.1142$. One multigrid V-cycle, with two
SOR pre-smoothing and two post-smoothing steps, is used to approximate the
background inverse. The rank tolerance is $10^{-10}$, and GMRES stops when the
relative residual is below $10^{-8}$, with at most 500 iterations. As stated
above, $B^{-1}$ is applied exactly in all tests.

Table~\ref{tab:reduction} shows the size of the reduced space and the accuracy of
the multigrid correction. The relative dimension falls from $18.1\%$ to
$11.8\%$, and both correction errors stay below $0.5\%$. The relative errors in Table~\ref{tab:reduction}(b) are computed by comparing the multigrid correction $C_1M_{\rm MG}C_1^T$ with the exact correction $C_1A_1^{-1}C_1^T$, in the Frobenius and spectral norms.

\begin{table}
\caption{Reduced-space size (left) and multigrid-correction error (right).}
\label{tab:reduction}
\centering
\begin{minipage}[t]{0.48\textwidth}
\centering\scriptsize
\textbf{(a) Dimensions}\par\smallskip
\renewcommand{\arraystretch}{1.00}
\begin{tabular}{|r|r|r|r|r|}
\hline
\multicolumn{1}{|c|}{$h$} &
\multicolumn{1}{c|}{$n_1$} &
\multicolumn{1}{c|}{$n_2$} &
\multicolumn{1}{c|}{$r$} &
\multicolumn{1}{c|}{$r/n_2$ (\%)}\\
\hline
0.8879 & 81   & 249   & 45   & 18.07\\
0.4522 & 289  & 977   & 149  & 15.25\\
0.2278 & 1089 & 3873  & 505  & 13.04\\
0.1142 & 4225 & 15425 & 1825 & 11.83\\
\hline
\end{tabular}
\end{minipage}\hfill
\begin{minipage}[t]{0.48\textwidth}
\centering\scriptsize
\textbf{(b) Relative error}\par\smallskip
\renewcommand{\arraystretch}{1.00}
\begin{tabular}{|r|r|r|}
\hline
\multicolumn{1}{|c|}{$h$} &
\multicolumn{1}{c|}{Frobenius} &
\multicolumn{1}{c|}{spectral}\\
\hline
0.8879 & $2.307\!\times\!10^{-4}$ & $2.763\!\times\!10^{-4}$\\
0.4522 & $2.334\!\times\!10^{-3}$ & $2.776\!\times\!10^{-3}$\\
0.2278 & $3.701\!\times\!10^{-3}$ & $3.704\!\times\!10^{-3}$\\
0.1142 & $4.417\!\times\!10^{-3}$ & $4.006\!\times\!10^{-3}$\\
\hline
\end{tabular}
\end{minipage}
\end{table}

Table~\ref{tab:spectral-gmres} reports the spectral diagnostics and GMRES
counts. The condition number of $SS_{\rm MG}^{-1}$ remains close to one, even
at the larger contrast. The $B^{-1}$ baseline worsens under refinement, while
the multigrid Woodbury preconditioner needs only two or three iterations. The
exact reduced Woodbury inverse is the reference case and takes one iteration on
every level. In the right-hand table, entries are ordered as
$B^{-1}$/MG-WB/Exact-WB.

\begin{table}
\caption{Spectral diagnostics (left) and GMRES iterations (right).}
\label{tab:spectral-gmres}
\centering
\begin{minipage}[t]{0.69\textwidth}
\centering\scriptsize
\textbf{(a) \(T=SS_{\rm MG}^{-1}\)}\par\smallskip
\setlength{\tabcolsep}{1.5pt}
\renewcommand{\arraystretch}{1.00}
\begin{tabular}{|c|r|r|r|r|r|r|}
\hline
& \multicolumn{3}{c|}{$1{:}10$}
& \multicolumn{3}{c|}{$1{:}100$}\\
\cline{2-7}
\multicolumn{1}{|c|}{$h$} &
\multicolumn{1}{c|}{$\kappa_2$} &
\multicolumn{1}{c|}{$\norm{I-T}_2$} &
\multicolumn{1}{c|}{$\rho(I-T)$} &
\multicolumn{1}{c|}{$\kappa_2$} &
\multicolumn{1}{c|}{$\norm{I-T}_2$} &
\multicolumn{1}{c|}{$\rho(I-T)$}\\
\hline
.8879 & 1.0011 & .0011 & $8.84\!\times\!10^{-5}$
      & 1.0038 & .0038 & $1.04\!\times\!10^{-4}$\\
.4522 & 1.0154 & .0152 & $1.22\!\times\!10^{-3}$
      & 1.0509 & .0497 & $1.57\!\times\!10^{-3}$\\
.2278 & 1.0305 & .0300 & $3.51\!\times\!10^{-3}$
      & 1.1210 & .1143 & $4.32\!\times\!10^{-3}$\\
.1142 & 1.0470 & .0459 & $5.32\!\times\!10^{-3}$
      & 1.1626 & .1510 & $6.18\!\times\!10^{-3}$\\
\hline
\end{tabular}
\end{minipage}\hfill
\begin{minipage}[t]{0.3\textwidth}
\centering\scriptsize
\textbf{(b) GMRES (B/MG/Ex)}\par\smallskip
\renewcommand{\arraystretch}{1.00}
\begin{tabular}{|c|r|r|}
\hline
\multicolumn{1}{|c|}{$h$} &
\multicolumn{1}{c|}{$1{:}10$} &
\multicolumn{1}{c|}{$1{:}100$}\\
\hline
.8879 & 25/2/1 & 31/2/1\\
.4522 & 31/3/1 & 70/3/1\\
.2278 & 39/3/1 & 88/3/1\\
.1142 & 45/3/1 & 115/3/1\\
\hline
\end{tabular}
\end{minipage}
\end{table}

Table~\ref{tab:truncation-cost} gives the truncation and cost results. On the
finest grid, retaining almost all of the correction energy does not guarantee
fast convergence: even small singular directions can matter to the
preconditioner. The growing LU fill ratio is another reason to treat the exact
correction only as a reference. Since the multigrid matrices were exported
before the MATLAB projection, these data are not an end-to-end timing study.

\begin{table}
\caption{Truncation experiment (left) and exact-background LU fill-in (right).}
\label{tab:truncation-cost}
\centering
\begin{minipage}[t]{0.49\textwidth}
\centering\scriptsize
\textbf{(a) Finest grid: \(n_2=15425\)}\par\smallskip
\renewcommand{\arraystretch}{1.10}
\begin{tabular}{|r|r|r|r|}
\hline
& & \multicolumn{2}{c|}{GMRES}\\
\cline{3-4}
\multicolumn{1}{|c|}{$k$} &
\multicolumn{1}{c|}{energy (\%)} &
\multicolumn{1}{c|}{$1{:}10$} &
\multicolumn{1}{c|}{$1{:}100$}\\
\hline
891  & 99.995718  & 41 & 105\\
1336 & 99.999363  & 29 & 78\\
1603 & 99.999958  & 24 & 60\\
1776 & 100.000000 & 3  & 3\\
1781 & 100.000000 & 3  & 3\\
\hline
\end{tabular}
\end{minipage}\hfill
\begin{minipage}[t]{0.49\textwidth}
\centering\scriptsize
\textbf{(b) Sparse LU for \(A_1\)}\par\smallskip
\renewcommand{\arraystretch}{1.10}
\begin{tabular}{|r|r|r|r|}
\hline
\multicolumn{1}{|c|}{$n_1$} &
\multicolumn{1}{c|}{$\operatorname{nnz}(A_1)$} &
\multicolumn{1}{c|}{$\operatorname{nnz}(L+U)$} &
\multicolumn{1}{c|}{fill}\\
\hline
81   & 393   & 868    & 2.21\\
289  & 1913  & 7436   & 3.89\\
1089 & 8409  & 62236  & 7.40\\
4225 & 35225 & 510524 & 14.49\\
\hline
\end{tabular}
\end{minipage}
\end{table}
The reduced multigrid correction captures the background response in the active
multiplier coupling space and produces behavior close to that of the exact
background solve. A scalable approximation of $B^{-1}$ remains the next step.
%======================================================================
% Bibliography
%======================================================================
% 
\bibliographystyle{splncs04}
\bibliography{References}
\end{document}